# Order restricted inference for comparing the cumulative incidence of a competing risk over several populations

**Hammou El Barmi**[*1], **Subhash Kochar**[2] and **Hari Mukerjee**[3]

*City University of New York, Portland State University and Wichita State University*

**Abstract:** There is a substantial literature on testing for the equality of the cumulative incidence functions associated with one specific cause in a competing risks setting across several populations against specific or all alternatives. In this paper we propose an asymptotically distribution-free test when the alternative is that the incidence functions are linearly ordered, but not equal. The motivation stems from the fact that in many examples such a linear ordering seems reasonable intuitively and is borne out generally from empirical observations. These tests are more powerful when the ordering is justified. We also provide estimators of the incidence functions under this ordering constraint, derive their asymptotic properties for statistical inference purposes, and show improvements over the unrestricted estimators when the order restriction holds.

## 1. Introduction

In the competing risks model, an experimental unit or subject is exposed to several risks, but the actual failure is attributed to exactly one of several possible distinct and exclusive types, called causes of failure. The available data is in the form of $(T, \delta)$, where $T$ is the time of failure (from any of the causes) and $\delta \in \{1, 2, \ldots\}$ is the cause of failure. Let $S$ denote the survival function (SF) of $T$. To determine the probability of failure by time $t$ from a certain cause, one defines the cumulative incidence function (CIF), a subdistribution function (SDF), by

$$F_i(t) = P[T \le t, \delta = i] = \int_0^t S(u) d\Lambda_i(u), \tag{1.1}$$

$i = 1, 2, \ldots$, where $\Lambda_i$ is the cumulative cause specific hazard function corresponding to the $i$th risk. We assume that $S$ is continuous.

Suppose that we have $k$ populations and in each population the units are exposed to several competing risks including one common risk whose severity we wish to compare across the populations. Without loss of generality, let the risk from cause 1 be of interest and let the remaining causes be lumped as other causes. For $i = 1, 2, \ldots, k$ and $j = 1, 2$, let $F_{ij}(t)$ be the cumulative incidence function due to cause $j$ in population $i$. Let $T_{il}$ be the time of failure of the $l$th unit in the $i$th population and let $\delta_{il}$ be the indicator of the corresponding cause of failure. Without censoring, the available data is in the form of independent random samples, $((T_{i1}, \delta_{i1}), (T_{i2}, \delta_{i2}), \ldots, (T_{in_i}, \delta_{in_i})), i = 1, 2, \ldots, k$. If $C_{il}$ is a censoring variable acting on the $l$th unit of population $i$ and $L_{il} = \min\{T_{il}, C_{il}\}$, then the data is in the form of $\{(L_{il}, \delta_{il})\}$, where $\delta_{il} \in \{0, 1, 2\}$ with $\delta_{il} = 0$ corresponding to censoring of this observation. For censored data, Gray [7] and Pepe and Mori [14] considered the test of

$$H_0 : F_{11}(t) = F_{21}(t) = \cdots = F_{k1}(t) \quad \text{for all } t \tag{1.2}$$

against all alternatives, using an integral of weighted differences of the empirical hazard rates and of the empirical CIF's, respectively. Note that it is not necessary for the other causes to be the same across the populations. Luo and Turnbull [13] also developed an integral test for equality of all CIF's across the populations when the competing risks are the same for each population. This paper also contains many other related references. Lin [11] developed a 2-sided Kolmogorov-Smirnov (KS) type test for the case $k = 2$ for the same problem and provided procedures for simultaneous confidence intervals for $F_{i1},\ i = 1, 2$. This test is asymptotically distribution-free and consistent against all alternatives. Since the asymptotic distribution of the test statistic is intractable, Lin [11] proposed a resampling technique for his test of (1.2) against all alternatives when $k = 2$.

*Supported by City University of New York through PSC-CUNY.
[1]Department of Statistics and CIS, Baruch College, City University of New York, New York, 10010, USA, e-mail: Hammou_ElBarmi@baruch.cuny.edu
[2]Department of Mathematics and Statistics, Portland State University, PO Box 751, Portland Oregon 97207-0751, USA, e-mail: kochar@pdx.edu
[3]Department of Mathematics and Statistics, Wichita State University, Wichita, Kansas 67260-0033, USA, e-mail: mukerjee@math.twsu.edu

In many of the examples given in the literature, intuitively we would expect that $\{F_{i1}\}$ is linearly ordered and the empirical evidence seems to support this. In this paper we consider a sequential KS type test of $H_0$ in (1.2) against $H_1 - H_0$, where

$$H_1 : F_{11} \le F_{21} \le \cdots \le F_{k1}. \tag{1.3}$$

This test is more powerful than the test against all alternatives when the linear ordering holds, as is to be expected. We also provide estimators of $\{F_{i1}\}$ under the restriction, derive their asymptotic properties, and provide improved simultaneous confidence band procedures under the order restriction.

For two linearly ordered distribution functions (DF's), Brunk et al. [2] found the nonparametric maximum likelihood estimator (NPMLE) whose asymptotic distribution was derived by Præstgaard and Huang [16]. This asymptotic distribution is too complicated for further statistical analyses. The NPMLE for the $k$-sample case for continuous DF's is unknown. Hogg [9] had suggested a simple estimator that can be applied to $k$ linearly ordered DF's. For the 2-sample case, Rojo and Ma [20] and Rojo [19] showed that Hogg's estimator generally has a lower mean square error (MSE) than the NPMLE at all quantiles. Rojo [18] and Rojo [19] showed that the limiting distributions of the estimators are the same as for the unrestricted case for $k = 2$ when the linear ordering is strict. El Barmi and Mukerjee [3] (hereafter referred to as EBM [3]) studied the limiting distributions for a general $k$, and when the ordering is not necessarily strict, providing procedures for asymptotic simultaneous confidence bands and KS type hypothesis tests. It turns out that all of their results hold when the DF's are replaced by SDF's in the uncensored case. For the censored case, however, the asymptotic distributions of the test statistics are intractable when the DF's are replaced by SDF's. In contrast, a similar KS type test statistic has an asymptotic distribution that is the supremum of a rescaled Brownian Motion when testing whether all the CIF's of a single population are equal (El

Barmi et al. [5], El Barmi and Mukerjee [4]). We use the method in Lin [11] for a sequential 1-sided KS type test in the $k$-sample case. The two main purposes of this paper are to extend Lin's test in Lin [11] for the 2-sample case to the $k$-sample case and to demonstrate how statistical inferences can be improved when the order restriction (1.3) holds.

In Sections 2 and 3 we consider the uncensored and the censored cases, respectively. In each case we provide estimators of $F_{11}, F_{21}, \ldots, F_{k1}$ that satisfy (1.3), prove their consistency, derive their asymptotic properties, show the improvements over the empiricals in terms of asymptotic MSE (AMSE), and provide improved confidence bands and hypothesis testing procedures. The properties of the estimators can be proven by simply replacing the DF in EBM [3] by an SDF. However, the confidence band and hypothesis testing procedures in the censored case require different treatments. In Section 4 we give an example of our procedures with real life data. In Section 5 we make some concluding remarks.

The symbols $\stackrel{w}{\Rightarrow}, \stackrel{d}{\rightarrow}$, and $\stackrel{d}{=}$ denote converges weakly to, converges in distribution to, and is equal in distribution to, respectively.

## 2. The uncensored case

### 2.1. The estimators

Consider independent random samples from $k$ life distributions, each subject being exposed to two causes of failure — cause 1, common to all $k$ populations, and cause 2, all "other causes," that need not be the same for all populations. Thus, the CIF's for the $i$th population are $F_{i1}$ and $F_{i2}$, given by (1.1). Let $T_{il}$ be the time of failure of the $l$th unit in the $i$th population and let $\delta_{il} \in \{1, 2\}$ be the corresponding cause of failure. Let $(T_{i1}, \delta_{i1}), (T_{i2}, \delta_{i2}), \ldots, (T_{in_i}, \delta_{in_i}), i = 1, 2, \ldots, k$, denote the samples.

The empirical estimate of $F_{ij}$ is given by

$$(2.1) \qquad \hat{F}_{ij}(t) = \frac{1}{n_i} \sum_{l=1}^{n_i} I(T_{il} \le t, \delta_{il} = j), \; j = 1, 2, \; i = 1, 2, \ldots, k.$$

Peterson [15] shows that this is also the NPMLE that is clearly strongly uniformly consistent. For a fixed $t$, let $\hat{\mathbf{F}}_1(t) = (\hat{F}_{11}(t), \hat{F}_{21}(t), \ldots, \hat{F}_{k1}(t))^T$ and let $\mathbf{n} = (n_1, n_2, \ldots, n_k)^T$. The restricted estimator of $F_{i1}(t)$ is obtained from the isotonic regression of $\hat{\mathbf{F}}_1$ with respect to the weight vector $\mathbf{n}$:

$$(2.2) \qquad \hat{F}^*_{i1}(t) = E_{\mathbf{n}}((\hat{F}_{11}, \ldots, \hat{F}_{k1})^T | \mathcal{A})_i(t), \quad 1 \le i \le k,$$

where $\mathcal{A} = \{\mathbf{u} \in R^k : u_1 \le u_2 \le \cdots \le u_k\}$, and $E_{\mathbf{w}}(\mathbf{u}|\mathcal{A})$ denotes the least squares projection of $\mathbf{u}$ onto $\mathcal{A}$ with the weight vector $\mathbf{w}$. It follows from Robertson et al. [17] (herafter referred to as RWD [17]) that

$$(2.3) \qquad \hat{F}^*_{i1} = \max_{r \le i} \min_{s \ge i} Av_{\mathbf{n}}[\hat{\mathbf{F}}_1; r, s]$$

where

$$Av_{\mathbf{n}}[\hat{\mathbf{F}}_1; r, s] = \sum_{l=r}^{s} n_l \hat{F}_{l1} / N_{rs}, \text{ where } N_{rs} = \sum_{l=r}^{s} n_l.$$

RWD [17] gives a comprehensive treatment of isotonic regression. It can be easily verified from the properties of isotonic regression that $\hat{F}^*_{i1}$ is nondecreasing, right continuous, and takes values in [0,1].

### 2.2. Asymptotic properties

From Corollary B, page 42 of RWD [17], we have

$$\max_{1\le i\le k} |\hat{F}^*_{i1}(t) - F_{i1}(t)| \le \max_{1\le i\le k} |\hat{F}_{i1}(t) - F_{i1}(t)| \quad \text{for each } t.$$

The following strong uniform consistency is an easy consequence of this.

**Theorem 2.1.** *If $F_{11}, F_{21}, \dots, F_{k1}$ satisfy (1.3) then*

$$P[||\hat{F}^*_{i1} - F_{i1}|| \to 0 \textit{ as } n_i \to \infty,\ i = 1, 2, \dots, k] = 1. \tag{2.4}$$

Next, we consider the asymptotic distributions of the restricted estimators. Define

$$Z_{i1n_i} = \sqrt{n_i}[\hat{F}_{i1} - F_{i1}] \quad \text{and} \quad Z^*_{i1\mathbf{n}} = \sqrt{n_i}[\hat{F}^*_{i1} - F_{i1}],\ i = 1, 2, \dots, k, \tag{2.5}$$

$$n = \sum_{i=1}^{k} n_i,\ \gamma_{in} = \frac{n_i}{n},\ \text{and assume that } \lim_{n\to\infty} \gamma_{in} = \gamma_i > 0 \text{ for all } i.$$

It follows from Breslow and Crowley [1] that

$$(Z_{11n_1}, Z_{21n_2}, \dots, Z_{k1n_k})^T \overset{w}{\Rightarrow} (Z_{11}, Z_{21}, \dots, Z_{k1})^T,$$

a $k$-variate Gaussian process with independent components, where $Z_{i1} \overset{d}{=} B_i(F_{i1})$ for all $i$. Here the $B_i$'s are independent standard Brownian Bridges. Let

$$\tilde{Z}_{i1n_i} = \sqrt{n}[\hat{F}_{i1} - F_{i1}] = Z_{i1n_i}/\sqrt{\gamma_{in}} \quad \text{and} \quad \tilde{Z}_{i1} = Z_{i1}/\sqrt{\gamma_i},\ i = 1, 2, \dots, k.$$

Then

$$(\tilde{Z}_{11n_1}, \tilde{Z}_{21n_2}, \dots, \tilde{Z}_{k1n_k})^T \overset{w}{\Rightarrow} (\tilde{Z}_{11}, \tilde{Z}_{21}, \dots, \tilde{Z}_{k1})^T. \tag{2.6}$$

We first consider the convergence in distribution at a fixed point, $t$. Define

$$\mathcal{S}_{it} = \{j : F_{j1}(t) = F_{i1}(t)\},\ i = 1, 2, \dots, k. \tag{2.7}$$

Note that $\mathcal{S}_{it}$ is a set of consecutive integers from $\{1, 2, \dots, k\}$ with $F_{j1}(t) - F_{i1}(t) = 0$ if $j \in S_{it}$, and, as $n \to \infty$,

$$\sqrt{n}[F_{j1}(x) - F_{i1}(t)] \to \infty,\ j > \mathcal{S}_{it},\ \text{and } \sqrt{n}[F_{j1}(t) - F_{i1}(t)] \to -\infty,\ j < \mathcal{S}_{it}, \tag{2.8}$$

where $j < (>)\mathcal{S}_{it}$ means $j < (>)l$ for all $l \in \mathcal{S}_{it}$. For $1 \le r \le s \le k$, let

$$\Gamma_{rsn} = \sum_{j=r}^{s} \gamma_{jn} = N_{rs}/n,\ \text{and let } \Gamma_{rs} = \lim_{n\to\infty} \Gamma_{rsn} = \sum_{j=r}^{s} \gamma_j.$$

**Theorem 2.2.** *Assume that $F_{11}(t) \le F_{21}(t) \le \cdots \le F_{k1}(t)$ for a fixed $t$ , then*

$$(Z^*_{11\mathbf{n}}(t), Z^*_{21\mathbf{n}}(t), \dots, Z^*_{k1\mathbf{n}}(t))^T \overset{d}{\to} (Z^*_{11}(x), Z^*_{21}(x), \dots, Z^*_{k1}(x))^T,$$

*where*

$$Z^*_{i1}(t) = \sqrt{\gamma_i} \max_{r\le i} \min_{i\le s} \frac{\sum_{\{r\le j\le s,\, r,s\in\mathcal{S}_{it}\}} \gamma_j \tilde{Z}_{j1}(t)}{\Gamma_{rs}}. \tag{2.9}$$

The idea of the proof comes from the fact that the $\{\mathcal{S}_{it} : 1 \le i \le k\}$ are disjoint "level sets" of $\{F_{i1}(t)\}$, and the isotonic estimator of $\hat{\mathbf{F}}_1$ in (2.3) will possibly violate the ordering only within each of these level sets asymptotically. The details are given in the proof of Theorem 2 of EBM [3].

To simplify the notation, for the rest of the paper we assume that the supports of $F_i$ is $[0, \tau_i]$ for some $0 < \tau_i < \infty$ (or $[0, \infty)$) for all $i$. As discussed in EBM after the proof of Theorem 2, weak convergence of the starred process on $\prod_i [0, \tau_i]$ requires some restrictions. We refer the reader to this paper for possible failure of convergence at isolated points without the following restrictions. Let

$$\mathcal{S}_i = \{j : F_{j1} = F_{i1}\} \text{ for } i = 1, 2, \ldots, k.$$

Consider the assumption

$$\inf_{\eta \le x \le \tau_i - \eta} [F_{j1}(x) - F_{i1}(t)] > 0 \text{ for all } \eta > 0 \text{ and } j > S_i, \quad i = 1, 2, \ldots, k-1. \tag{2.10}$$

The extension of point-wise convergence in distribution to weak convergence using assumption (2.10) is straightforward (see proof of Theorem 4 in EBM [3]).

**Theorem 2.3.** *If (1.3) and (2.10) hold. Then*

$$(Z^*_{11\mathbf{n}}, Z^*_{21\mathbf{n}}, \ldots, Z^*_{k1\mathbf{n}})^T \stackrel{w}{\Rightarrow} (Z^*_{11}, Z^*_{21}, \ldots, Z^*_{k1})^T,$$

*where*

$$Z^*_{i1} = \sqrt{\gamma_i} \max_{r \le i} \min_{i \le s} \frac{\sum_{\{r \le j \le s,\, r,s \in \mathcal{S}_i\}} \gamma_j \tilde{Z}_{j1}}{\Gamma_{rs}}.$$

Note that, if $S_i = \{i\}$, then $Z^*_{i1\mathbf{n}} \stackrel{w}{\Rightarrow} Z_{i1}$ under the conditions of Theorem 2.3 as should be the case.

### *2.3. Comparison with empirical estimators*

In this section we compare some of the properties of our estimators with those of the unrestricted empiricals.

From Kelly's analysis in Kelly [10], we have the following interesting result. Suppose that $\mathbf{X} = (X_1, X_2, \ldots, X_l)^T$ has independent components, $X_i \sim N(0, \sigma_i^2)$, $\sigma_i^2 > 0,\ i = 1, 2 \ldots, l$. Define the weight vector $\mathbf{w}$ by

$$\mathbf{w} = \left(\frac{1}{\sigma_1^2}, \frac{1}{\sigma_2^2}, \ldots, \frac{1}{\sigma_l^2}\right)^T,$$

and let $\mathbf{X}^* = E_{\mathbf{w}}(\mathbf{X}|\mathcal{A})$ be the isotonic regression of $\{X_i\}$ with weights $\{1/\sigma_i^2\}$, constrained to lie in $\mathcal{A}$ as in (2.2). Then

$$P(|X_i^*| \le u) > P(|X_i| \le u) \text{ for all } u > 0,\ i = 1, 2, \ldots, l. \tag{2.11}$$

Now suppose that $\mathcal{S}_{it}$, as defined in (2.7), has more than one element for some $t$ with $0 < F_{i1}(t) < 1$. Then, from the form of (2.9) in Theorem 2.2, $\{Z^*_{j1}(t)/\sqrt{\gamma_j} : j \in S_{it}\}$ is the least squares projection of $\{\tilde{Z}_{j1}(t) : j \in S_{it}\}$ with (equal) weights $\{\sigma_i^{-2}(t) : j \in S_{it}\}$. Since $\tilde{Z}_j(t) \sim N(0, F_{i1}(t)[1 - F_{i1}(t)]/\gamma_i)$, we have the following theorem.

**Theorem 2.4.** *Under the conditions of Theorem 2.3, if $S_{it}$ contains more than one element for some $t$ with $0 < F_{i1}(t) < 1$, then*

$$P[|Z^*_{j1}(t)| \leq u] > P[|Z_{j1}(t)| \leq u] \tag{2.12}$$

*for all $u > 0,\ j \in S_{it}$.*

Two immediate consequences of this theorem are that $E[Z^*_{i1}(t)]^2 < E[Z_{i1}(t)]^2$, implying an improvement in AMSE, and a possible increase of confidence coefficients for confidence bands centered around the restricted estimators as opposed to the empiricals with the same bandwidths. These results can be substantially sharpened when $k = 2$. We refer the reader to Section 4.2 in EBM [3].

### 2.4. Hypothesis testing

Here we consider testing $H_0$ against $H_1 - H_0$, where $H_0$ and $H_1$ are given by (1.2) and (1.3), respectively. For $k = 2$ we could use the 1-sided KS test. For the $k$-sample test, Hogg [9] suggested a sequential pairwise test as follows. We test $H_{0j}$ against $H_{1j} - H_{0j}$, where, for $j = 2, \ldots, k-1$,

$$\begin{aligned} H_{0j} &: F_{11} = F_{21} = \cdots = F_{j1} \text{ and} \\ H_{1j} &: F_{11} = F_{21} = \cdots = F_{j-1,1} \leq F_{j1}. \end{aligned} \tag{2.13}$$

Let

$$\Gamma_{jn} = \sum_{i=1}^{j} \gamma_{in},\ 1 \leq j \leq k,\ \Delta_{jn} = \gamma_{jn}\Gamma_{j-1,n}/\Gamma_{jn},\ j \geq 2,$$

The test statistic, $\mathcal{T}_n$, is defined as

$$\mathcal{T}_n = \max_{2\leq j\leq k} \sup_t \mathcal{T}_{jn}(t),\ \text{ where } \mathcal{T}_{jn} = \sqrt{n}\sqrt{\Delta_{jn}}[\hat{F}_{j1} - Av_{\mathbf{n}}[\hat{\mathbf{F}}_{\mathbf{1}}; 1, j-1]].$$

The statistic $\sup_t \mathcal{T}_{jn}(t)$ is used to test $H_{0j}$ against $H_{1j} - H_{01}$ for $2 \leq j \leq k$; $\mathcal{T}_n$ is the maximum of these quantities, and the test rejects $H_0$ for large values of $\mathcal{T}_n$.

Under $H_0$, denoting the common CIF by $F_{11}$, we have

$$\mathcal{T}_{jn} = \sqrt{\Delta_{jn}}[\tilde{Z}_{j1\mathbf{n}} - Av_{\gamma_{\mathbf{n}}}[\tilde{\mathbf{Z}}_{1\mathbf{n}}; 1, j-1]],$$

where $\tilde{\mathbf{Z}}_{1\mathbf{n}} = (\tilde{Z}_{11\mathbf{n}}, \tilde{Z}_{21\mathbf{n}}, \ldots, \tilde{Z}_{k1\mathbf{n}})^T$ and $\gamma_n = (\gamma_{1n}, \gamma_{2n}, \ldots, \gamma_{kn})^T$. By the continuous mapping theorem, $(\mathcal{T}_{2n}, \mathcal{T}_{3n}, \ldots, \mathcal{T}_{kn})^T \stackrel{w}{\Rightarrow} (\mathcal{T}_2, \mathcal{T}_3, \ldots, \mathcal{T}_k)^T$, where

$$\mathcal{T}_j = \sqrt{\Delta_j}[\tilde{Z}_{1j} - Av_{\boldsymbol{\gamma}}[\tilde{\mathbf{Z}}_1; 1, j-1]],$$

$\boldsymbol{\gamma} = (\gamma_1, \gamma_2, \ldots, \gamma_k)^T$, $\tilde{\mathbf{Z}}_1 = (\tilde{Z}_{11}, \tilde{Z}_{21}, \ldots, \tilde{Z}_{k1})^T$, with $\tilde{Z}_{i1}$ as defined in (2.6) , and $\Delta_j = \lim_n \Delta_{jn}$, that is equal to $\gamma_j\Gamma_{j-1}/\Gamma_j$ with $\Gamma_j = \sum_{i=1}^{j} \gamma_i$. By computing covariances, it can be shown that the $\mathcal{T}_j$'s are independent, and that

$$\mathcal{T}_j \stackrel{d}{=} B_j(F_{11}), \quad 2 \leq j \leq k.$$

Therefore $\mathcal{T}_n$ converges in distribution to $\mathcal{T}$, where

$$\mathcal{T} = \max_{2\leq j\leq k} \sup_t \mathcal{T}_j(t).$$

and

$$\begin{aligned} P(\mathcal{T} \geq t) &= 1 - P(\sup_t \mathcal{T}_j(t) < t, j = 2, \ldots, k) \\ &= 1 - \prod_{j=2}^{k} P(\sup_t B_j(F_{11}(t)) < x) \\ &\leq 1 - (1 - e^{-2t^2})^{k-1} \approx (k-1)e^{-2t^2}, \end{aligned}$$

the Bonferroni approximation without the independence of $\{\mathcal{T}_{jn}\}$ which is quite good unless $k$ is fairly large.

## 3. The censored case

In this section we consider the same set-up as in Section 2 except that there may be censoring in addition to the competing risks. We assume that the censoring mechanisms are independent of the competing risks. We denote the causes of failure by 0, 1 and 2, where $\delta = 0$ corresponds to the observation being censored. We assume that the censoring variable corresponding to $T_{ij}$ is $C_{ij}$ with a continuous SF $G_i$ for all $i$ and $j$. In this setting, we observe $\{(L_{ij}, \delta_{ij})^T\}$, where $L_{ij} = T_{ij} \wedge C_{ij}$ and $\delta_{ij} = \delta_{ij} I(T_{ij} \leq C_{ij})$. Therefore $\{L_{ij},\ j = 1, 2, \ldots, n_i\}$ is a random sample from $H_i = 1 - S_i G_i,\ i = 1, 2, \ldots, k$. Let $\pi_i = S_i G_i$ and let $\tau_m = \min_i \tau_i$, where $\tau_i$ is the right endpoint of the support of $F_{i1}$ as defined before.

### *3.1. Estimators and consistency*

In order to use the martingale formulation for asymptotics, we define the unconstrained estimator of $F_{i1}$ by

$$\hat{F}_{i1}(t) = \int_0^t \hat{S}_i(u) d\hat{\Lambda}_{i1}(u) \tag{3.1}$$

where, for technical reasons, $\hat{S}_i$ is the left continuous Kaplan-Meier estimator of $S_i$ and $\hat{\Lambda}_{ij}$ is the Nelson-Aalen estimator of $\Lambda_{ij}$, given by (Fleming and Harrington [6])

$$\hat{\Lambda}_{ij}(t) = \sum_{l=1}^{n_i} \frac{I[L_{il} \leq t, \delta_{il} = j]}{\sum_{k=1}^{n_i} I[L_{ik} \geq L_{il}]}.$$

As in the uncensored case, our restricted estimator of $F_{i1}(t)$ is given by

$$\hat{F}^*_{i1}(t) = E_{\mathbf{n}}((\hat{F}_{11}, \ldots, \hat{F}_{k1})^T | \mathcal{A})_i(t), \quad 1 \leq i \leq k, \tag{3.2}$$

where $\mathbf{n}$ and $\mathcal{A}$ are as defined in Section 2. It is well known that $(||h||_0^b = \sup_{0 \leq t \leq b} |h(t)|)$

$$P[||\hat{F}_{i1} - F_{i1}||_0^b \to 0 \text{ as } n_i \to \infty, i = 1, 2, \ldots, k] = 1 \tag{3.3}$$

for any $b < \tau_m$. Thus, as in the uncensored case, we have

$$P(||\hat{F}^*_{i1} - F_{i1}||_0^b \to 0 \text{ as } n_i \to \infty, i = 1, 2, \ldots, k) = 1 \tag{3.4}$$

from the properties of isotonic regression.

### 3.2. Weak convergence

Let $Z_{i1n_i} = \sqrt{n_i}[\hat{F}_{i1} - F_{i1}]$ and $Z^*_{i\mathbf{n}} = \sqrt{n_i}[\hat{F}^*_{i1} - F_{i1}],\ 1 \le i \le k$. From Lin [11],

$$(Z_{11n_i}, Z_{21n_2}, \ldots, Z_{k1n_k})^T \overset{w}{\Rightarrow} (Z_{11}, Z_{21}, \ldots, Z_{k1})^T \quad \text{on } \prod_i [0, \tau_i), \tag{3.5}$$

where $\{Z_{i1}\}$ are independent mean-zero Gaussian processes with the covariances given by (for $s \le t$)

$$\begin{aligned} Cov(Z_{i1}(s), Z_{i1}(t)) = & \int_0^s [1 - F_{i1}(s) - F_{i2}(u)][1 - F_{i1}(t) - F_{i2}(u)] \frac{d\Lambda_{i1}(u)}{\pi_i(u)} \\ & + \int_0^s [F_{i1}(s) - F_{i1}(u)][F_{i1}(t) - F_{i1}(u)] \frac{d\Lambda_{i2}(u)}{\pi_i(u)}. \end{aligned}$$

By independence, $Cov(Z_{i1}(s), Z_{j1}(t)) = 0$ for $i \ne j$.

### 3.3. Improvements by restricted estimation

Kelly's result in Kelly [10] leading to Theorem 2.4 depends only on normality and independence of $\{Z_{i1}\}$. Thus, Theorem 2.4 holds in the censored case also. This implies that the AMSE of the restricted estimators could be reduced by using the restricted estimators and the confidence coefficients in asymptotic confidence bands centered around the restricted estimators could be more than those centered around the unrestricted estimators using the same bandwidths.

### 3.4. Confidence bands and hypothesis testing

The distributions of $Z_{i1}$'s are intractable. To approximate, Lin [11], based on an earlier paper by Lin et al. [12], develops the following martingale formulation. In the following, the indices $i$, $j$ and $l$ range from 1 to $k$, 1 to $n_i$, and 1 to 2, respectively. Let $Y_{ij}(t) = I(L_{ij} \ge t)$, $N_{ijl}(t) = I(L_{ij} \le t, \delta_{ij} = l)$, $M_{ijl}(t) = N_{ij}(t) - \int_0^t Y_{ij}(u) d\Lambda_{il}(u)$, $\bar{Y}_i(t) = \sum_{j=1}^{n_i} Y_{ij}(t)$, $\bar{N}_{il}(t) = \sum_{j=1}^{n_i} N_{ijl}(t)$, and $\bar{M}_{il}(t) = \sum_{j=1}^{n_i} M_{ijl}(t)$. In the counting process notation, the data $\{(L_{ij}, \delta_{ij})\}$ is represented by $\{(Y_{ij}(\cdot), N_{ijl}(\cdot))\}$. Lin [11] shows that, by replacing $\bar{M}_{il}$ by $\sum_{j=1}^{n_i} V_{ijl} N_{ijl}$, where $V_{ijl}$'s are independent standard normals, the conditional distribution of

$$\begin{aligned} \hat{Z}_{i1}(t) \equiv \sqrt{n}_i \sum_{j=1}^{n_i} \Bigg[ & \int_0^t \frac{[1 - \hat{F}_{i2}(u)] V_{ijl} dN_{ijl}(u)}{\bar{Y}_i(u)} + \int_0^t \frac{\hat{F}_{ij} V_{ij2} dN_{ij2}(u)}{\bar{Y}_i(u)} \\ & - \bar{F}_{i1}(t) \int_0^t \frac{V_{ijl} dN_{ij1}(u) + V_{ij2} dN_{ij2}(u)}{\bar{Y}_i(u)} \Bigg] \end{aligned} \tag{3.6}$$

given $\{(Y_{ij}(\cdot), N_{ijl}(\cdot))\}$ is asymptotically equivalent to the unconditional distribution of $Z_{i1}(\cdot)$. This clever observation allows one to simulate probabilities involving $Z_{i1}(t)$ by generating random samples of $V_{ijl}$ while fixing $\{(Y_{ij}(\cdot), N_{ijl}(\cdot))\}$ at their observed values. Lin [11] considered various forms of confidence bands using the transformed processes

$$D_{i1}(t) = \sqrt{n}_i g(t) [\phi(\hat{F}_{i1}(t)) - \phi(F_{i1}(t))],$$

where $\phi$ is a known function with a continuous non-zero derivative $\phi'$ and $g$ is a weight function. By the functional delta method theorem, $D_{i1}(t)$ is asymptotically equivalent to $g(t)\phi'(F_{i1}(t))Z_{i1}(t)$, that can be approximated by $\hat{D}_{i1}(t) = g(t)\phi'(\hat{F}_{i1}(t))\hat{Z}_{i1}(t)$. If $P[\sup_{t_1 \le t \le t_2} |\hat{D}_{ij}(t)| > q_\alpha] = \alpha$, then

$$\phi(\hat{F}_{i1}(t)) \pm n^{-1/2} q_\alpha / g(t),\ t_1 \le t \le t_2$$

is a $(1-\alpha)$ confidence band for $\phi(F_{i1}(\cdot))$ on $[t_1, t_2]$. This confidence band procedure can be carried out for each population individually. As noted above in Section 3.3, one could possibly improve on the coverage probabilities with the same bandwidths if the bands are centered at the restricted estimators.

One could carry out the hypothesis testing in Section 2.4 with censored observations in exactly the same way using the same notation, but using the estimators given by (3.1) instead of (2.1). The limit of the process $\mathcal{T}_{jn}$ used in testing $H_{0j}$ against $H_{1j} - H_{0j}$ will still be given by

$$(3.7)\quad \mathcal{T}_{jn} = \sqrt{\Delta_{jn}}[\tilde{Z}_{j1\mathbf{n}} - Av_{\gamma_\mathbf{n}}[\tilde{\mathbf{Z}}_{1n}; 1, j-1]] \stackrel{w}{\Rightarrow} \mathcal{T}_j = \sqrt{\Delta_j}[\tilde{Z}_{1j} - Av_{\boldsymbol{\gamma}}[\tilde{\mathbf{Z}}_1; 1, j-1]]$$

under $H_{0j}$. The distributions of the $\mathcal{T}_j$'s are intractable and they are generally not independent in the censored case. However, we can use the resampling scheme using (3.6) to approximate the distribution of $\tilde{Z}_{1j}$ by $\hat{Z}_{j1}/\sqrt{\gamma_{jn}}$ to carry out the pairwise tests and use the Bonferroni bound given at the end of Section 2.4 to compute the $p$-value when $k$ is not large, say less than or equal to 10.

## 4. An example

To illustrate the theoretical results, we analyze a set of mortality data provided by Dr. H. E. Walburg, Jr., of the Oak Ridge National Laboratory and reported by Hoel [8] The data were obtained from a laboratory experiment on 82 RFM strain male mice who were kept in a germ-free environment (group 1) and another 99 who were kept in a conventional laboratory environment (group 2). Each group received a radiation dose of 300 rads at 5-6 weeks of age. After autopsy, the causes of death were classified as thymic lymphoma (cause 1), reticulum cell sarcoma (cause 2), and other causes which we consider here to be a censoring mechanism. Hoel [8] writes, "Biologists believe that both of these two diseases are lethal and that they are independent of one another and of other causes of death." Here, we need not assume the independence of the causes of death, but we do assume that they are independent of other causes. Intuitively, we would expect that the CIFs associated with all risks to be smaller in the germ-free environment. Since cause 1 is considered to be most the lethal of all the causes, we test $H_0 : F_{11} = F_{21}$ against $H_1 - H_0$, where $H_1 : F_{11} \le F_{21}$, using the test statistic described in Section 3.4. The test statistic we use is $\sup_t \mathcal{T}_{2n}(t)$, where $\mathcal{T}_{jn}$ is as defined in (3.7) . Since the limiting distribution of this test statistic is not tractable, we use the technique described in Section 3.4 to simulate 10,000 replicates from the distribution $\mathcal{T}_{jn}$. The value of the test statistic in this case 1.11 corresponding to a $p$-value of 0.676. This does not provide any evidence for the ordering assumption, possibly because thymic lymphoma is extremely lethal, and the germ-free environment does not provide enough protection against it. However, for estimation purposes the ordering still seems reasonable. The empirical and the order restricted estimators are displayed in Figures 1 and 2, respectively. In this case, the restricted estimators are exactly the same.

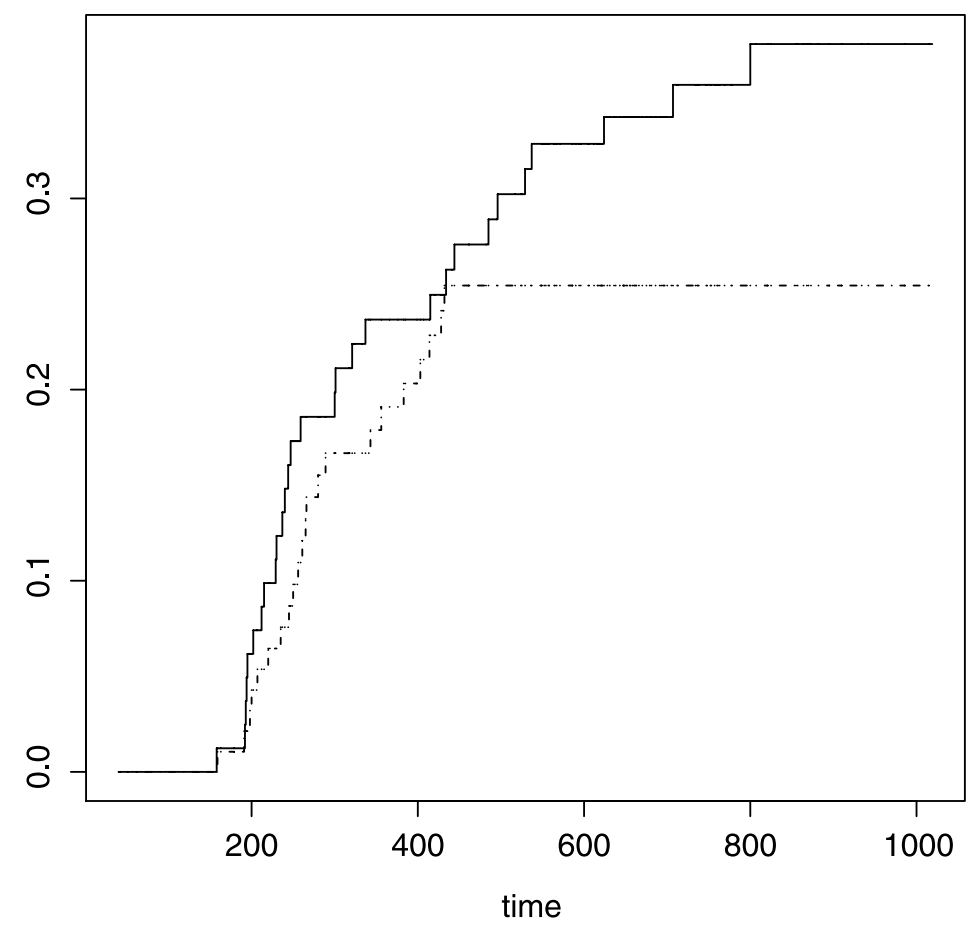

FIG 1. *Unrestricted estimators of the cumulative incidence functions* $\hat{F}_{11}$ *(dotted line) and* $\hat{F}_{21}$ *(solid line).*

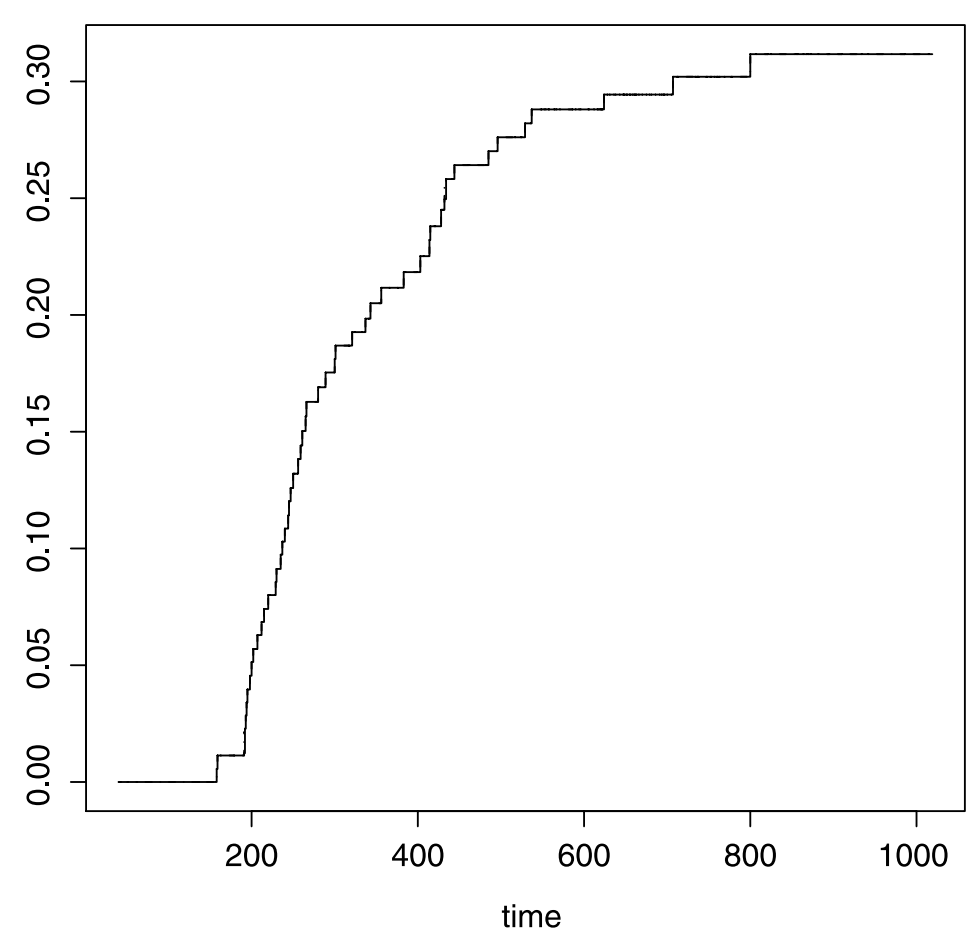

FIG 2. *Restricted estimators of the cumulative incidence functions* $\hat{F}^*_{11}$ *(dotted line) and* $\hat{F}^*_{21}$ *(solid line).*

## 5. Concluding remarks

In competing risks studies, there are several tests in the literature for equality of the cumulative incidence function due to a specific cause across several populations against all alternatives. Most of these are integral tests. Lin [11] proposed a 2-sided Kolmogorov-Smirnov type test in the case of two populations and against all alternatives. In some situations it is reasonable to assume that the incidence functions of interest in the various populations are linearly ordered. For the 2-sample case, the precedure in Lin [11] could be easily modified to a 1-sided test. This is more powerful when the ordering holds. We show how to extend this to the $k$-sample case using a sequence of 1-sided tests. We also provide estimators of the incidence functions under the order restriction, derive their asymptotic properties, and show improvements in the inference procedures over the unrestricted estimators when the order restriction is obeyed, and even when it is violated by small amounts.

We give a real life example to illustrate our procedure.

**Acknowledgments.** The authors are grateful to the Editor and an anonymous referee for valuable comments that have led to an improved presentation.

## References

[1] Breslow, N. and Crowley, J. (1974). A large sample study of the life table and product limit estimates under random censorship. *Ann. Statist.* **2** 437–453. MR0458674

[2] Brunk, H. D., Franck, W. E., Hanson, D. L. and Hogg, R. V. (1966). Maximum likelihood estimation of the distributions of two stochastically ordered random variables. *J. Amer. Statist. Assoc.* **61** 1067–1080. MR0205367

[3] El Barmi, H. and Mukerjee, H. (2005). Inferences under a stochastic ordering constraint: The $k$-sample case. *J. Amer. Statist. Assoc.* **100** 252–261. MR2156835

[4] El Barmi, H. and Mukerjee, H. (2006). Restricted estimation of the cumulative incidence functions corresponding to competing risks. In *Lecture Notes Monogr. Ser.* **49** 241–252. Institute of Mathematical Statistics. MR2338546

[5] El Barmi, H., Kochar, S. C., Mukerjee, H. and Samaniego, F. J. (2004). Estimation of cumulative incidence functions in competing risks studies under an order restriction. *J. Statist. Plann. Inference* **118** 145–165. MR2015226

[6] Fleming, T. R. and Harrington, D. P. (1991). *Counting Processes and Survival Analysis*. Wiley, New York. MR1100924

[7] Gray, R. J. (1988). A class of $K$-sample tests for comparing the cumulative incidence of a competing risk. *Ann. Statist.* **16** 1141–1154. MR0959192

[8] Hoel, D. G. (1972). A representation of mortality data by competing risks. *Biometrics* **28** 475–478.

[9] Hogg, R. V. (1962). Iterated tests of the equality of several distributions. *J. Amer. Statist. Assoc.* **61** 579–585. MR0159373

[10] Kelly, R. E. (1989). Stochastic reduction of loss in estimating normal means by isotonic regression. *Ann. Statist.* **17** 937–940. MR0994278

[11] Lin, D. Y. (1997). Non-parametric inference for cumulative incidence functions in competing risks studies. *Statistics in Medicine* **16** 901–910.

[12] Lin, D. Y., Fleming, T. R. and Wei, L. J. (1994). Confidence bands for survival curves under the proportional hazards model. *Biometrika* **81** 73–81. MR1279657

[13] Luo, X. and Turnbull, B. W. (1999). Comparing two treatments with multiple competing risks endpoints. *Statist. Sinica* **9** 985–997.

[14] Pepe, M. and Mori, M. (1993). Kaplan-Meier, marginal or conditional probability curves in summarizing competing risks failure time data. *Statistics in Medicine* **12** 737–751.

[15] Peterson, Jr., A. V. (1977). Expressing the Kaplan-Meier estimator as a function of empirical subsurvival functions. *J. Amer. Statist. Assoc.* **72** 854–858. MR0471165

[16] Præstgaard, J. T. and Huang, J. (1996). Asymptotic theory for nonparametric estimation of survival curves under order restrictions. *Ann. Statist.* **24** 1679–1716. MR1416656

[17] Robertson, T., Wright, F. T. and Dykstra, R. L. (1988). *Order Restricted Statistical Inference*. Wiley, New York. MR0961262

[18] Rojo, J. (1995). On the weak convergence of certain estimators of stochastically ordered survival functions. *J. Nonparametr. Statist.* **4** 349–363. MR1366781

[19] Rojo, J. (2004). On the estimation of survival functions under a stochastic order constraint. In *The First Erich L. Lehmann Symposium—Optimality. IMS Lecture Notes Monogr. Ser.* **44** 37–61. Inst. Math. Statist., Beachwood, OH. MR2118560

[20] Rojo, J. and Ma, Z. (1996). On the estimation of stochastically ordered survival functions. *J. Statist. Comput. Simulation* **55** 1–21. MR1700892